\documentclass{amsart}
\usepackage{graphicx, amssymb, amsmath, amsthm}

\newcommand {\C} {\mathbb{C}}

\newcommand {\Hom} {\mathrm{Hom}}
\newcommand {\End} {\mathrm{End}}
\newcommand {\Tr}{\mathrm{Tr}}

\newtheorem{thm}{Theorem}[section]
\newtheorem{lemma}[thm]{Lemma}

\newtheorem{defn}[thm]{Definition}

\begin{document}
\subjclass[2000]{}
\title{Mednykh's Formula via Lattice Topological Quantum Field Theories}
\author{Noah Snyder}
\email{nsnyder@math.berkeley.edu}
\address{Department of Mathematics \\ University of California at Berkeley \\ Berkeley, California 94720}

\begin{abstract}
Mednykh \cite{Mednykh} proved that for any finite group $G$ and any orientable surface $S$, there is a formula for $\# \Hom(\pi_1(S), G)$ in terms of the Euler characteristic of $S$ and the dimensions of the irreducible representations of $G$.  A similar formula in the nonorientable case was proved by Frobenius and Schur \cite{F-S}.  Both of these proofs use character theory and an explicit presentation for $\pi_1$.  These results have been reproven using quantum field theory (\cite{FQ}, \cite{M-Y}, and others).  Here we present a greatly simplified proof of these results which uses only elementary topology and combinatorics.  The main tool is an elementary invariant of surfaces attached to a semisimple algebra called a lattice topological quantum field theory.
\end{abstract}

\maketitle
\section{Introduction}
A lattice topological quantum field theory is a topological invariant of surfaces attached to a semisimple algebra (more generally, a knowledgeable Frobenius algebra \cite{LP}) which is computed using an explicit triangulation.  These invariants are called topological quantum field theories because they behave nicely under gluing \cite{Atiyah}.  Lattice topological quantum field theories were originally introduced in \cite{FHK} as a toy model for understanding Turaev-Viro invariants of $3$-manifolds \cite{TV}.  Although topological invariants of surfaces are uninteresting on their own (since the Euler characteristic is so successful), there has been a resurgence of interest in $2$-dimensional topological quantum field theories due to their appearance in Khovanov homology \cite{Khovanov}.  Lattice topological quantum field theories are of particular interest in Khovanov homology because they extend easily to surfaces with corners \cite{LP2} and to unoriented surfaces \cite{TT}.  This paper gives another application of these invariants: a rapid and elementary proof of Mednykh's formula.

To describe Mednykh's formula, let's fix some notation.  Let $G$ be a finite group and $S$ be a closed surface.  Let $\chi$ denote the Euler characteristic and let $\hat{G}$ be the set of isomorphism classes of irreducible representations of $G$.  Let $\nu$ be the Frobenius-Schur indicator.  (For the definition of the Frobenius-Schur indicator and group theory background see Section 3.)  Let $d(V)$ denote the dimension of $V$.

The goal of this paper is to prove the following formulas from \cite{Mednykh} and \cite{F-S} (see Section $2$ for more history).
If $S$ is orientable, then \begin{equation} \label{orientedresult} \sum_{V \in \hat{G}} d(V)^{\chi(S)} = \#G^{\chi(S)-1} \# \Hom(\pi_1(S), G).\end{equation}
If $S$ is non-orientable, then \begin{equation}\label{unorientedresult} \sum_{V \in \hat{G}} (\nu(V)d(V))^{\chi(S)} = \#G^{\chi(S)-1} \# \Hom(\pi_1(S), G).\end{equation}

The main advantage of our approach is that the proof (and the formulas themselves) can be easily reconstructed from a single sentence:
\begin{quote}\emph{``Compute the lattice topological quantum field theory invariant of $S$ attached to the group algebra with respect to the two obvious bases."}\end{quote}

The first basis is the group-like elements and yields the right hand sides of the main equations.  The second basis is the matrix elements of the irreducible representations and yields the left hand sides of the main equations.  We address the oriented and non-oriented cases seperately, because in the latter case the definition of a lattice TQFT requires extra data.

The plan of the paper is as follows.  Section $2$ explains the history of Mednykh's formula and its many proofs.  Section $3$ is a refresher on  the necessary notions from elementary group theory.  Section $4$ gives the definition of a lattice topological quantum field theory attached to a semisimple algebra.  The main differences from \cite{FHK} is that we emphasize the basis-free nature of this invariant and we generalize it to nonorientable surfaces.  This generalization requires an involutive $*$-structure on the semisimple algebra.  In Section $5$ we explicitly compute this invariant with respect to the two different bases of the group algebra.

I would like to thank Chris Schommer-Pries and Vladmir Turaev for several very helpful conversations.

\section{The History of Mednykh's Formula}

Formulas \ref{orientedresult} and \ref{unorientedresult} have been rediscovered and reproven many times over the years.  Since most papers seem to be unaware of the full history, we have done our best to collect this history here.  The reader uninterested in history may easily skip to the next section.

The original argument in the non-orientable  case, due to \cite{F-S}, is purely algebraic.  It uses generators and relations for $\pi_1$ and character theory techniques.  According to \cite[p. 51]{Cavlieri}  Burnside used similar techniques to prove the orientable case of Mednykh's formula in the special case of the symmetric group.  These algebraic techniques are sufficient for the general orientable case as was proved by \cite{Mednykh} (see also \cite[\S 4.]{M-Y2} and \cite[Lemma 4.2.4]{subfactor} for details).  This proof is rapid and elementary, but its reliance on generators and relations for $\pi_1$ obscures the relationship with topology.

Mednykh's formula was rediscovered in the early 90s from a quantum field theoretic perspective. The quantum field theoretic proofs come in two main flavors: (2+1)-dimensional and (1+1)-dimensional.  In the former case, one computes the dimension of the vector space associated with the boundary of a $3$-manifold.  In the latter case one computes the scalar associated to a closed surface.

The main sources for the (2+1)-dimensional approach are \cite{Verlinde}, \cite{DW}, and \cite{FQ}.  The mathematically inclined reader will find the last reference easier going.  These proofs compare a gauge-theoretic computation with a gluing based computation.  From the (2+1)-dimensional perspective, one may think of Mednykh's formula as a special case of Verlinde's formula applied to the quantum double of the group ring.  See \cite[Chapter 4.2]{subfactor} for the relationship between Verlinde's formula and Mednykh's formula.

See \cite[\S 1.2]{Segal} for a sketch of a (1+1)-dimensional approach inspired by \cite{DW} which also compares a gauge theoretic computation with a gluing computation (or see \cite[p. 78]{Bartlett} for further details).  A completely different $2$-dimensional approach is taken in \cite{M-Y} and \cite{M-Y2}, where they use matrix integrals.  Of the field theoretic techniques, to our knowledge only the Mulase-Yu approach has been adapted to the non-orientable case (although, see \cite{AN} for some results in this direction).

From our perspective, each of the above field theoretic proofs is unnecessarily long and complicated.  By contrast, our proof does not use any serious geometry, does not require a $3$-dimensional invariant (nor the Hopf algebra theory required for $3$-dimensions), requires no familiarity with physics, and only uses one construction of the invariant.  We hope that this simplicity will allow these beautiful formulas (and their relation to quantum topology) to be understood by a wider audience.

\section{Group Theory Background}

Let $G$ be a finite group.  Let $k$ be an algebraically closed field of characteristic relatively prime to $\#G$.  (The reader will lose very little by assuming that $k = \C$ since a corollary of Brauer's theorem states that the dimensions of the simple $G$-modules is independent of the choice of $k$.)  Consider the group algebra $k[G]$.  Let $\hat{G}$ denote a set of representatives of the isomorphism classes of irreducible representations of $G$.  By Mashke's theorem $k[G]$ is semisimple, so by Artin-Wedderburn, we have that (using $M_n$ to denote $n$-by-$n$ matrices) \begin{equation} \label{A-W} k[G] \cong \bigoplus_{V \in \hat{G}_k} M_{d(V)}.\end{equation}  This decomposition defines a basis of matrix elements $e_{ij}(V)$ for $k[G]$.

The reader only interested in the orientable case may skip the rest of this section.  For the non-orientable case we will need to understand how the linear anti-involution $*$ on $k[G]$ defined by $g^* = g^{-1}$ (extended by linearity) acts on the right hand side of Equation \ref{A-W}.

\begin{defn}
If $A$ is a semisimple algebra with an anti-involution $*$, and $V$ is any finite-dimensional $A$-module.  Define $V^*$ to be the $A$-module defined by $af(v) = f(a^*v)$.
\end{defn}

If $A = k[G]$ and $*$ is defined as above, then this definition of $V^*$ agrees with the usual definition.  However, since the twisted group case is also of interest (see \cite{Turaev}), we give the definition in more generality.  Because $*$ is invertible, if $V$ is simple then so is $V^*$.  Since $*$ is an involution, we see that $*$ yields an involution on the set of isomorphism classes of simple $A$-modules.  Thus $*$ also gives an involution on the matrix factors $\End(V)$ of $A$.  We want to understand how this involution acts.

\begin{defn}
If $A$ is a semisimple algebra with an anti-involution $*$, and $V$ is any finite-dimensional $A$-module, then we define the Frobenius-Schur indicator $\nu(V)$ as follows.  If $V \ncong V^*$ then $\nu(V) = 0$.  If $V \cong V^*$ fix an isomorphism $f: V \rightarrow V^*$.  Consider $f^*: V = V^{**} \rightarrow V$.  We define the Frobenius-Schur indicator by $f^* = \nu(V) f$.
\end{defn}

\begin{defn}
For $i = \pm 1$, let $\hat{G}_i$ denote the subset of $\hat{G}$ of irreducible representations with $\nu(V) = i$.  Let $\hat{G}_0$ denote a choice of representative from each pair $V$, $V^*$ when $\nu(V) = 0$.
\end{defn}

Since $f=f^{**} = \nu(V)^2 f$ we see that $\nu(V) = \pm 1$.  If $\nu(V) = 0$, then $*$ interchanges $\End(V)$ and $\End(V^*)$.  We use  $(M_n \oplus M_n)^\text{swap}$ to denote this $*$-algebra.  If $\nu(V) = 1$, then up to change of basis $*$ is given by transpose.  We use $M_n$ to denote the $*$-algebra where $*$ is given by transpose.  If $\nu(V) = -1$ then up to change of basis $*$ is given by transpose and then conjugation by an anti-symmetric invertible matrix.  In particular, if $\nu(V) = -1$, then $d(V)$ is even.  We use $M^{\text{anti}}_{2n}$ to denote this antisymmetric $*$-structure.

Finally, we note that as $*$-algebras (where we define $(a\otimes b)^* = a^* \otimes b^*$),
$$(M_n \oplus M_n)^\text{swap} \cong M_n \otimes (k \oplus k)^\text{swap}$$
$$M^\text{anti}_{2n} \cong M_n \otimes M^{\text{ant}i}_2.$$

(In the specific case where $A = k[G]$ with the $*$ structure $g^* = g^{-1}$ and when the characteristic of $k$ is not $2$, there's an explicit formula for $\nu$ due to Frobenius and Schur saying that $\nu(V) = \sum_{g \in G} \chi_V(g^2).$  This formula is crucial for the character theoretic proofs of Mednykh's formula, but we will not use it.)

To summarize this section, we can refine the algebra isomorphism in Equation \ref{A-W} to the following isomorphism of $*$-algebras.

\begin{align} \label{starA-W}
& k[G] \cong \\
\nonumber \cong & \left(\bigoplus_{V \in \hat{G}_1} M_{d(V)}\right) \oplus \left(\bigoplus_{V \in \hat{G}_0}(M_{d(V)} \oplus M_{d(V)})^\text{swap} \right) \oplus \left(\bigoplus_{V \in \hat{G}_{-1}} M^\text{anti}_{d(V)}\right) \\
 \nonumber \cong & \left(\bigoplus_{V \in \hat{G}_1} M_{d(V)}\right) \oplus \left(\bigoplus_{V \in \hat{G}_0}M_{d(V)} \otimes (k \oplus k)^\text{swap} \right) \oplus \left(\bigoplus_{V \in \hat{G}_{-1}} M_{\frac{d(V)} {2}} \otimes M^\text{anti}_2\right)
\end{align}

\section{Topological Invariants and Semisimple Algebras}

In this section we define an invariant of surfaces attached to a semi-simple algebra.  This invariant can be extended to a $2$-dimensional topological field theory with corners (see \cite{LP} in the oriented case, and \cite{KM} in the unoriented case), but we will concentrate on the case of closed surfaces.  This invariant is called a lattice topological quantum field theory.  The definition here is largely identical to that in \cite{FHK}, but here we emphasize the basis-free nature of the invariant.  (See \cite{Baez} for the right motivation for these definitions.)

Let $\Tr(x): A \rightarrow k$ be the trace of $m_x$, the multiplication by $x$ map from $A$ to $A$.  Notice that the map $T_k: A^{\otimes k} \rightarrow k$ given by $\Tr(x_1 x_2 \cdots x_k)$ is invariant under cyclic permutations.  By semisimplicity, $T_2$ is a nondegenerate symmetric billinear form $A\otimes A \rightarrow k$.  This gives an identification $A \rightarrow A^*$.  Since $A$ is semisimple this map is invertible giving a map $A^* \rightarrow A$ which gives a symmetric map $p: k \rightarrow A \otimes A$.  We will use Sweedler's notation, $p(1) = \sum p_1 \otimes p_2$.  

\begin{defn}
A flag in a triangulated surface is a pair (face, edge) such that the edge is is contained in the face.
\end{defn}

\begin{defn}
Let $S$ be an oriented surface with a fixed triangulation.  Let $\#V$, $\#E$, $\#F$ be the numbers of vertices, edges, and faces respectively.  To each flag we associate a copy of $A$.  To each edge we assign the map $p: k \rightarrow A^{\otimes 2}$.  To every oriented triangle we associate the map $T_3: A^{\otimes 3}\rightarrow k$ (since this map is invariant under cyclic permutations this map only depends on the choice of orientation).  Thus to the triangulated surface we have assigned a map $$I_A: k = k^{\otimes \# E} \rightarrow A^{\otimes \text{flags}} \rightarrow k^{\otimes \# F} = k.$$
\end{defn}

\begin{thm}
$I_A(S)$ depends only on the topology of $S$ and is independent of the triangulation.
\end{thm}
\begin{proof}
Since this has already been proved in \cite{FHK} and \cite{LP} we just sketch the proof.

All equivalences of oriented triangulated surfaces are generated by the following two Pachner moves \cite{Pachner}, called the $1-3$ and $2-2$ moves.

\begin{center}
\includegraphics[height=.5in]{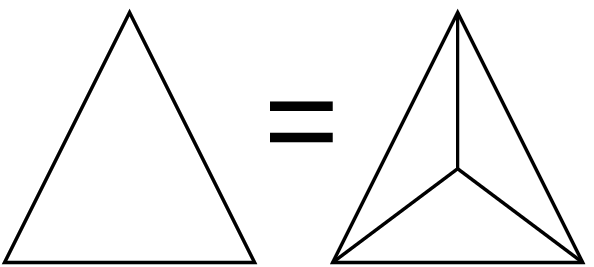}
\end{center}

\begin{center}
\includegraphics[height=.5in]{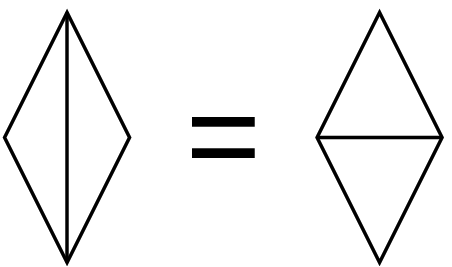}
\end{center}

Thus, in order to show that $I_A$ is a topological invariant, it is enough to show that the maps associated to those diagrams are equal.  The $2-2$ move reduces easily to associativity of multiplication.  The $1-3$ move follows from our choice of normalization for the trace ($\Tr(1) = d(A)$).
\end{proof}

The unoriented case is only slightly more complicated.  Now we consider a semisimple $*$-algebra.  That is we fix an involution $*$ such that $(ab)^* = b^*a^*$, $\Tr(a^*) = \Tr(a)$, and $\sum p_1 \otimes p_2^* = \sum p_1^* \otimes p_2$.  Take an unoriented triangulated surface together with a choice of orientation of each triangle.  Again we assign a vector space to each flag.  To each triangle we assign the map $T_3$ with the inputs ordered by the orientation on the triangle.  To each edge we assign $p$ if the orientations of the two triangles agree and $(1 \otimes *) \circ p$ if the triangles have opposite orientations.  Again we define $I_A$ to be the composition of all of these maps.

To see that $I_A$ does not depend on the choice of orientation of the individual triangles, we show that it does not change when you reverse the orientation on one triangle.  This move reverses the order of product in the inputs of $T_3$ and switches whether the orientations agree or disagree for each of the edges.  Thus, the equation $$\Tr(abc) = \Tr((abc)^*) = \Tr(c^*b^*a^*)$$ implies that $I_A$ is independent of the choice of orientations.

We collect a few elementary results about the invariant $I$.  If $A$ and $B$ are semisimple ($*$-)algebras then $A \otimes B$ and $A \oplus B$ are also semisimple ($*$-)algebras. An explicit computation shows that $I_{A \otimes B}(S) = I_A(S) I_B(S)$, that $I_{A \oplus B}(S) = I_A(S) + I_B(S)$, and that $I_k(S) = 1$.

\section{Computing the Invariant attached to a Group Algebra}

Since the characteristic of $k$ is relatively prime to $\#G$, by Maschke's theorem, $k[G]$ is a semisimple algebra.  Furthermore, there is a star structure given by $g^* = g^{-1}$.  Thus, the construction from the previous section yields an invariant of surfaces attached to a group algebra.  We compute this invariant with respect to two different bases.

\begin{thm}
For any surface (orientable or not), \[I_{k[G]}(S) = \#G^{\chi(S)-1} \#\Hom(\pi_1(S),G).\]
\end{thm}
\begin{proof}
This result is well-known among experts.  In $3$-dimensions it appears as early as \cite[p. 3]{Kuperberg}.  See \cite[Chapter 4.1]{subfactor} for a proof in $3$-dimensions that is very similar to our argument.

We compute $I_{k[G]}$ with respect to the basis of group-like elements to get a state-sum formula for $I_{k[G]}$.  In this case the computation is identical for the oriented and unoriented versions, so we do both computations together.  The map associated to each edge sends $1 \mapsto 1/\#G \sum g \otimes g^{\pm 1}$ where the sign is $-$ if the orientations agree and $+$ if they disagree.  The map $T_3$ sends $a \otimes b \otimes c$ to $\#G$ if $abc=1$ and $0$ otherwise.

Thus, $$I_{k[G]}(S) = \#G^{\#F - \#E} Z(G, S),$$ where $Z(G,S)$ is the number of ways of labeling each flag with an element of $G$ such that the two labels at each edge are equal if the orientation reverses and inverse if it does not reverse, and such that the product around each triangle is $1$.  Two adjacent triangles with the same orientation assign opposite orientations to their common edge.  Thus, equivalently one can label \emph{oriented} edges of $S$ with elements of $G$ such that the two orientations of an edge are assigned inverse group elements and the oriented product around every triangle is $1$.  We call such labelings consistent, and note that $Z(G,S)$ is the number of consistent labelings.

Fix $v_0$ a base vertex of $S$, and an oriented path $P_v$ along edges of $S$ from $v_0$ to $v$ for every other vertex $v$.  We construct a bijection between consistent labelings of $S$ and the set $G^{V-\{v_0\}} \times  \Hom(\pi_1(S), G)$ as follows.  Let $f$ be a consistent labeling thought of as a map from oriented edges to $G$.  To every vertex $v \neq v_0$ we assign the element $\prod_{e \in P_v} f(e)$, and to any loop $L$  we assign the element $f(L) = \prod_{e \in L} f(e)$.  Notice that this assignment to $L$ only depends on the class of $L$ in $\pi_1(S)$, because of the consistency condition on triangles.

Conversely if we have an element of $g$ assigned to each of the $P_v$ and each of the $L_i$ we can recover the consistent labeling, thus proving that this assignment is a bijection.  Consider the edge $E$ from $v$ to $v'$.  Let $L$ be the loop $P_{v'}^{-1} \circ E \circ P_{v}$.  We can recover $f(E) =f(P_v') f(L)f(P_v)^{-1}$.  See the figure for a picture of this bijection.

\begin{center}
\includegraphics[height=1.5in]{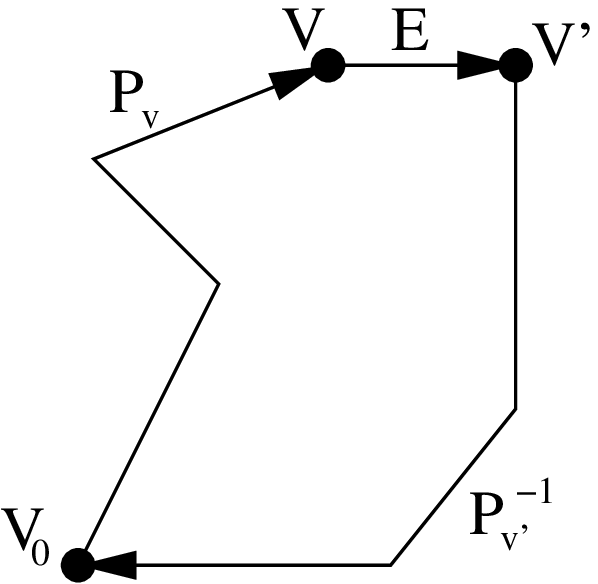}
\end{center}

Hence, we have computed that 
\begin{align}  \label{orientedrhs} I_{k[G]}(S) &= \#G^{\#F - \#E} Z(G,S) \\
\nonumber &= \#G^{\#F -\#E + \#V-1} \#\Hom(\pi_1(S), G) \\
\nonumber &= \#G^{\chi(S)-1} \#\Hom(\pi_1(S), G).\end{align}
\end{proof}

Now we would like to compute $I_{k[G]}$ using the decomposition as a direct sum of matrix algebras from Section $2$.  In order to do that we compute $I_{M_n}$.  

\begin{thm}
For any surface $S$ (orientable or not), \[I_{M_n}(S) = n^{\chi(S)}.\]
\end{thm}
\begin{proof}
The map associated to each edge sends $1$ to $1/n \sum_{i,j} e_{ij} \otimes e_{ji}$ if the orientations agree and   $1/n \sum_{i,j} e_{ij} \otimes e_{ij}$ if they do not agree.  The map $T_3$ sends $e_{ij} \otimes e_{jk} \otimes e_{ki}$ to $n$ and all triples not of that form to $0$.  Just as before we can think of this as labeling directed edges instead of flags.  Thus, $$I_{M_n}(S) = n^{\#F-\#E} Z(M_n, S),$$ where $Z(M_n, S)$ is the number of ways of labeling oriented edges of $S$ by pairs $(i,j)$ such that the same pair with opposite orientation is labeled $(j,i)$, and such that adjacent edges in the same oriented triangle are labeled $(i, j)$ and $(j, k)$.

Such labelings are clearly equivalent to labeling each vertex with a number, and giving each edge the numbers attached to its initial vertex and final vertex as in the figure.

\begin{center}
\includegraphics[height=1in]{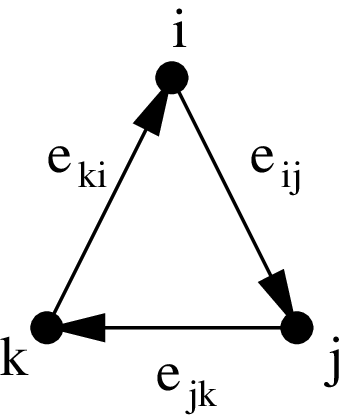}
\end{center}

Thus we see that $Z(M_n, S) = n^{\#V}$, and $I_{M_n} = n^{\chi(S)}$.
\end{proof}

If $S$ is orientable, then use the algebra isomorphism $k[G] \cong \bigoplus_V M_{d(V)}$ to get that \begin{equation} \label{orientedlhs} I_{k[G]} = \sum_V I_{M_{d(V)}} = \sum_V (d(V))^{\chi(S)}.\end{equation}  Combining Equation \ref{orientedlhs} with Equation \ref{orientedrhs} yields a proof of Equation \ref{orientedresult}.

If $S$ is not necessarily orientable, then we need to split up the calculation for the matrix algebras based on the $*$-structure.  Recall from Section $3$ that we have the isomorphism of $*$-algebras:
$$k[G] \cong \left(\bigoplus_{V \in \hat{G}_1} M_{d(V)}\right) \oplus \left(\bigoplus_{V \in \hat{G}_0}M_{d(V)} \otimes (k \oplus k)^\text{swap} \right) \oplus \left(\bigoplus_{V \in \hat{G}_{-1}} M_{\frac{d(V)} {2}} \otimes M^\text{anti}_2\right).$$

Thus all that remains to compute is $I_{(k\oplus k)^\text{swap}}(S)$ and $I_{M^\text{anti}_2}(S)$.  The former computation is completely elementary, but the latter is a bit trickier and uses significantly more topology than any other computation in this paper. 

\begin{thm} 
$I_{(k\oplus k)^\text{swap}}(S)$ counts the number of orientations of $S$.  In particular it is either $2$ or $0$ depending on whether $S$ is orientable.
\end{thm}
\begin{proof}
Let $e_1$ and $e_{-1}$ denote the two standard basis vectors of $k \oplus k$.  By definition $e_i^* = e_{-i}$.  The map associated to each edge sends $1$ to $e_1 \otimes e_1 + e_{-1} \otimes e_{-1}$ if the orientations agree and to $e_1 \otimes e_{-1} + e_{-1} \otimes e_{1}$ if the orientation disagrees.  The map associated to each triangle sends $e_i \otimes e_j \otimes e_k$ to $1$ if $i=j=k$ and to $0$ otherwise.  Thus, $I_{(k\oplus k)^\text{swap}}(S)$ counts the number of ways you can color each triangle with $e_{\pm 1}$ so that if two adjacent triangles have the same color when their orientations agree and the opposite colors when orientations differ.  Given any such coloring, one may reverse the orientation on all the triangles labeled with $e_{-1}$ to get an orientation of the surface.  Thus, $I_{(k\oplus k)^\text{swap}}(S)$ counts the number of orientations of $S$.
\end{proof}

\begin{thm} \label{quaternion}
$I_{M^\text{anti}_2}(S) = (-2)^{\chi(S)}$.
\end{thm}

In order to prove this theorem we will use the following lemma in the topology of surfaces.

\begin{lemma}
Let $S$ be a surface.  There exists a finite maximal collection of $n$ non-intersecting Mobius strips in $S$.  Furthermore, $n \equiv \chi(S) \mod 2$.
\end{lemma}
\begin{proof}
Suppose that $S$ is orientable, then $S$ contains zero Mobius strips.  Thus we need only prove that its Euler characteristic is even.  This follows from Poincar\'e duality.  Namely, there is a symplectic intersection form on $H_1(S, \C)$.  Thus $\dim H_1(S, \C)$ is even.  Hence, $\chi(S) = \dim H_2 - \dim H_1 + \dim H_0 = 2-\dim H_1$ is also even.

Now consider a non-orientable surface $S$.  Since it is non-orientable it contains a Mobius strip.  Removing this Mobius strip and gluing in a disc increases the Euler characteristic by $1$.  Since for any surface $S'$ we have that $$\chi(S') = \dim H_2 - \dim H_1 + \dim H_0 = 2-\dim H_1 \leq 2,$$ we may repeat this process at most $2-\chi(S)$ times.  Hence there exists a finite maximal collection, whose size we call $n$, of non-intersecting Mobius strips in $S$.  Removing all of these Mobius strips and replacing them with discs results in a surface $S''$ which is oriented (because a Mobius strip in $S''$ would contradict maximality).  Since replacing a Mobius strip with a disc increases Euler characteristic by $1$, we see that $\chi(S) = \chi(S'') -n$.  By the first pagraph $\chi(S'')$ is even, hence, $n \equiv \chi(S) \mod 2$.
\end{proof}

\begin{proof}[Proof of Theorem \ref{quaternion}]
First we work out what $*$ does on an explicit basis of $M^\text{anti}_2$, namely

$$\left(\begin{array}{cc} a & b \\ c & d \end{array}\right)^* = \left(\begin{array}{cc} 0 & 1 \\ -1 & 0 \end{array}\right) \left(\begin{array}{cc} a & b \\ c & d \end{array}\right)^T \left( \begin{array}{cc} 0 & 1 \\ -1 & 0 \end{array}\right)^{-1} = \left(\begin{array}{cc} d & -b \\ -c & a \end{array}\right)
$$

We'll refer to the two standard basis vectors as $v_1$ and $v_{-1}$, thus the matrix elements are $e_{\pm 1, \pm 1}$.  With respect to this $*$-structure, $e_{ij}^* = i j e_{-j, -i}$.

The map attached to each edge sends $1$ to $\frac{1}{2} \sum e_{ij} \otimes e_{ji}$ if the orientations agree, and to $\frac{1}{2} \sum i j e_{i j} \otimes e_{-i, -j}$ if they disagree.  The map attached to each triangle sends $e_{ij} \otimes e_{jk} \otimes e_{ki}$ to $2$, and all other basis vectors to $0$.

As before with respect to this basis we see that $2^{\#E - \#F} I_{M^\text{anti}_2}(S)$ counts certain colorings of the surface, but this time we need to count with a sign.  Define a corner to be a pair (face, vertex) such that the face contains the vertex.  We call two corners adjacent if they have the same vertex, and the two faces share a common edge.  Then, $2^{\#E - \#F} I_{M^\text{anti}_2}(S)$ counts the number of labelings of each corner with $\pm 1$, such that if two corners are adjacent, then the labels agree if and only if the orientations of their faces agree.  These labelings are counted with a sign which is $(-1)^m$ where $m$ is the number of edges of $S$ such that the two faces containing that edge have opposite orientation, and the two corners of the same face adjacent to that edge have opposite labels.

The following picture shows part of a consistent labeling of the corners, together with the corresponding assignments of basis vectors to flags.  In this picture only the bottom edge contributes a sign (the left edge has the same orientation on both sides, while the right edge has the two corners on the same face having the same label). 

\begin{center}
\includegraphics[height=3.5in]{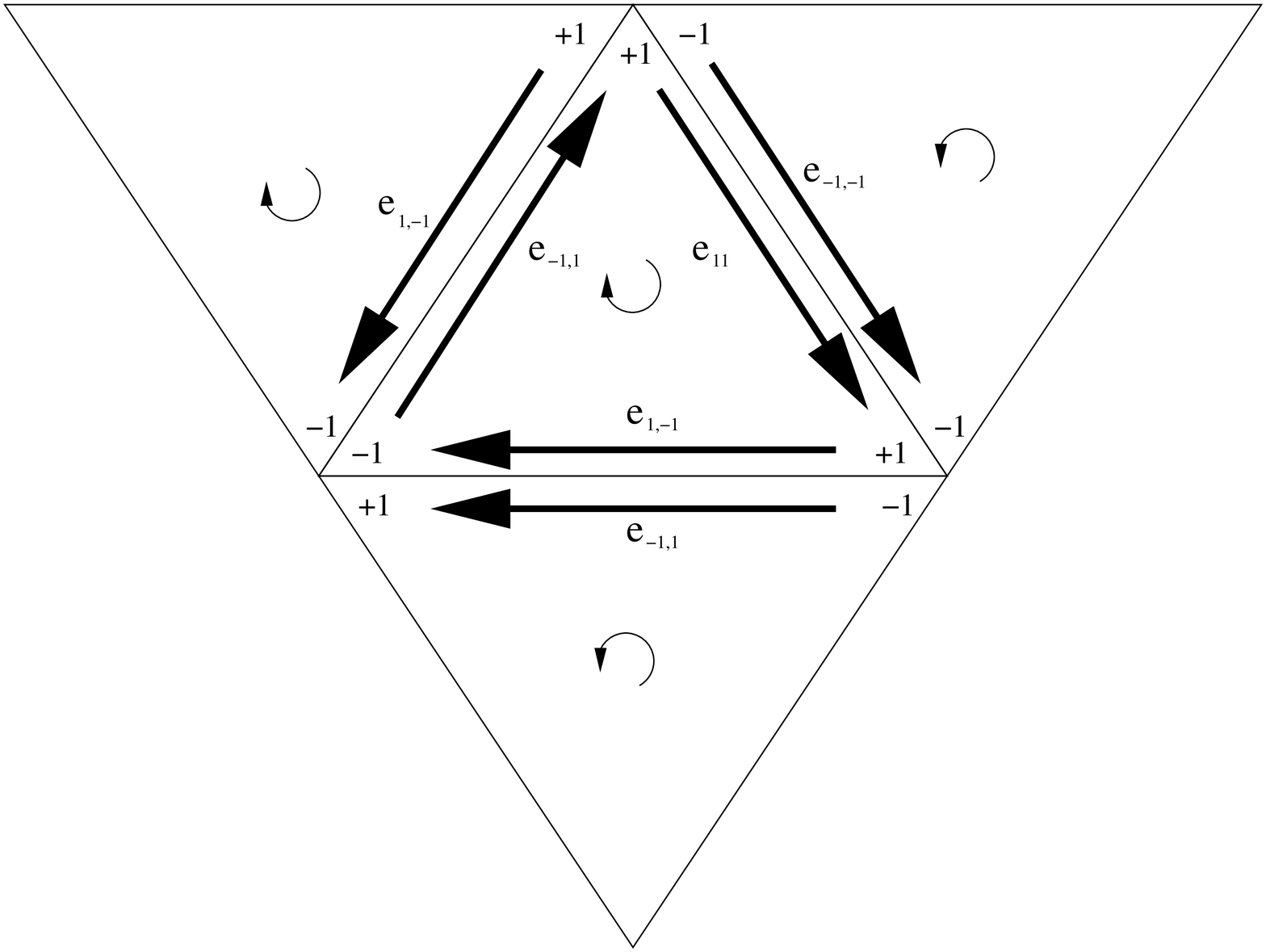}
\end{center}

If we fix a coloring of one corner at every vertex, then there is exactly one way to extend this to a compatable labeling.  Thus the total number of labels is $2^{\# V}$.  Furthermore, if you change the coloring of one of these chosen corners, it is easy to see that this does not change the sign (because there are an even number of orientation changes going around each vertex).  Hence, $I_{M^\text{anti}_2}(S) = \pm 2^{\chi(S)}$.  We need only determine which sign to use for which surfaces.

Consider a maximal collection of $n$ nonintersecting Mobius strips in $S$.  It is easy to see that after possibly barycentric-subdividing $n$-times, we can assume that these Mobius strips are tubular neighborhoods of non-intersecting edge loops.   Since removing all these Mobius strips results in an orientable surface with boundary, we may choose orientations on the triangles such that the orientation only changes when it crosses one of these $n$ loops.  

The following figure shows a typical tubular neighborhood of one of these $n$ loops with a consistant coloring.  Only the edge marked with an asterisk contributes a sign.  Hence we see that each loop contributes a $-1$.

\begin{center}
\includegraphics[height=1.25in]{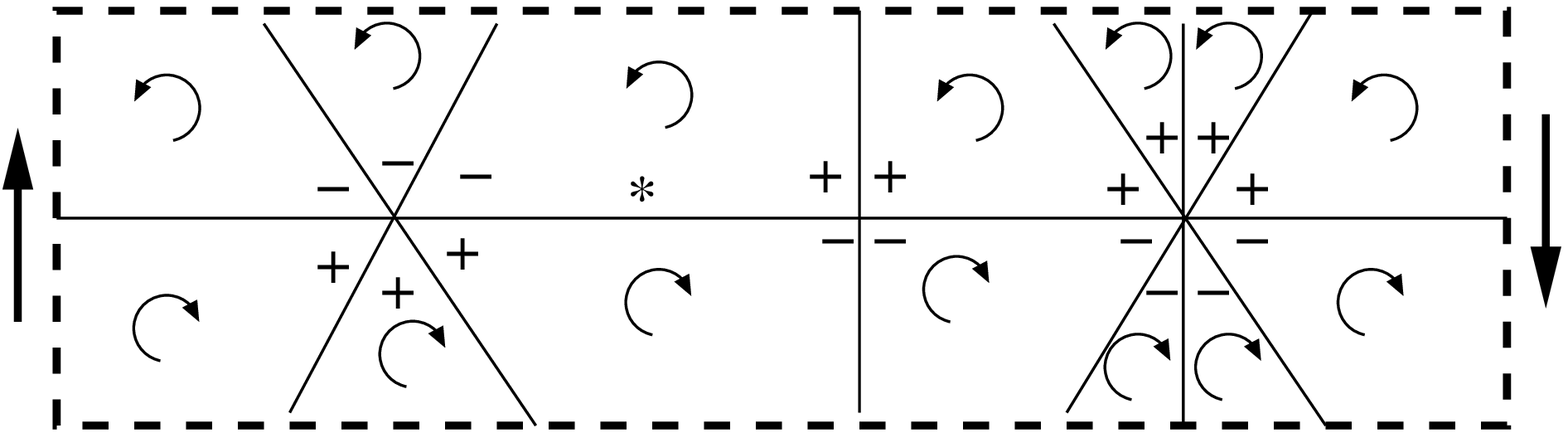}
\end{center}

Thus, the sign of any coloring of $S$ is $(-1)^n$ where $n$ is the size of a maximal collection of disjoint Mobius strips.  Using the lemma, we conclude that $$I_{M^\text{anti}_2}(S) = 2^{\#F-\#E+\#V} (-1)^{n} = (-2)^{\chi(S)}.$$
\end{proof}

In conclusion, if $S$ is unorientable we see that
\begin{align}
\label{unorientedlhs} & I_{k[G]}(S) = \\
\nonumber &= \sum_{V \in \hat{G}_1} I_{M_{d(V)}}(S) + \sum_{V \in \hat{G}_0} I_{M_{d(V)}}(S) I_{(k \oplus k)^\text{swap}}(S) + \sum_{V \in \hat{G}_{-1}} I_{M_{\frac{d(V)} {2}}}(S) I_{M^\text{anti}_2}(S) \\
\nonumber &= \sum_{V \in \hat{G}_1} d(V)^{\chi(S)} + \sum_{V \in \hat{G}_0} (d(V))^{\chi(S)}\cdot 0 + \sum_{V \in \hat{G}_{-1}} (d(V)/2)^{\chi(S)} (-2)^{\chi(S)} \\
\nonumber &= \sum_{V \in \hat{G}_1} d(V)^{\chi(S)}  + \sum_{V \in \hat{G}_{-1}} (-d(V))^{\chi(S)} = \sum_{V \in \hat{G}} (\nu(V) d(V))^{\chi(S)}. \end{align}

Combining Equation \ref{unorientedlhs} with Equation \ref{orientedrhs} yields a proof of Equation \ref{unorientedresult}.

\end{document}